\documentstyle[12pt]{article}

\font\tenmsb=msbm10 scaled \magstep1
\font\sevenmsb=msbm7 scaled \magstep1
\font\fivemsb=msbm5 scaled \magstep1

\newfam\msbfam
\textfont\msbfam=\tenmsb
\scriptfont\msbfam=\sevenmsb
\scriptscriptfont\msbfam=\fivemsb
\def\Bbb#1{{\fam\msbfam\relax#1}}

\begin{document}

\newtheorem{th}{Theorem}[section]
\newtheorem{cor}{Corollary}[section]
\newtheorem{prop}{Proposition}[section]

\def\hR{\hat{\Bbb R}^m}
\def\hG{\hat\Gamma}
\def\R{{\Bbb R}}
\def\C{{\Bbb C}}
\def\F{{\Bbb F}}
\def\N{{\Bbb N}}
\def\Z{{\Bbb Z}}
\def\F{{\Bbb F}}
\def\T{{\Bbb T}}
\def\CP{{\Bbb CP}^k}
\def\be{\begin{equation}}
\def\en{\end{equation}}
\def\Subset{\subset\subset}
\def\<{\langle }
\def\>{\rangle }
\def\O{\Omega}
\def\o{\omega}
\def\G{\Gamma}
\def\g{\gamma}
\def\e{\varepsilon}
\def\d{\delta}
\def\wz{\wedge_\Z}
\def\p{\varphi}
\def\l{\lambda}
\def\a{\alpha}
\def\bs{~\hfill\rule{7pt}{7pt}}
\def\ws{$\Box$ }
\def\fxy{\bar f(\tau\cdot jx,y)}
\def\pxy{\p(\tau\cdot jx,y)}
\def\supp{{\rm supp}}
\def\const{{\rm const}}
\def\sp{{\rm sp\,}}
\def\Arg{{\rm Arg }}
\def\dd{(i/2\pi)\partial\bar\partial}
\def\Lin{{\rm Lin}}
\def\fHT{f\in APH(T_\Omega)}
\def\fAT{f\in AP(T_\Omega)}
\def\AM{APM_{1,1}(T_\O)}
\def\hRO{\hR\times\O}
\def\L{\Lambda}
\def\f{\bar f(\tau,y)}
\def\grad{{\rm grad}}
\def\bp{\bar\partial}
\def\dist{{\rm dist}}
%%%%%%%%%%%%%%%%%%%%%%%%%%%%%%%%%%%%%%%%%%%%%%%%%%%%%%%%%%%%%%%%
\title{
    \bf Almost periodic divisors, holomorphic functions, and holomorphic
mappings \footnote {This research was supported by INTAS-99-00089
project.} }

\author{Favorov S.Ju.}

\date{}

\maketitle

\begin{abstract}

We prove that to each almost periodic, in the sense of
distributions, divisor $d$ in a tube $T_\O\subset\C^m$ one can
assign a cohomology class from $H^2(K,\Z)$ (actually, the first
Chern class of a special line bundle over $K$ generated by
$d$) such that the trivial cohomology class represents the
divisors of all almost periodic holomorphic functions on $T_\O$;
here $K$ is the Bohr compactification of $\R^m$. This description
yields various geometric conditions for an almost periodic
divisor to be the divisor of a holomorphic almost periodic
function. We also give a complete description for the divisors of
homogeneous coordinates for holomorphic almost periodic curves; in
particular, we obtain a description for the divisors of
meromorphic almost periodic functions.

\medskip


\medskip


\end{abstract}

\medskip
\section*{}\label{sec0}
\bigskip

The classical theory of almost periodic functions has found a lot
of applications in various branches of mathematics, from
differential equations (cf. \cite{So}, \cite{Fi}, \cite{P}) to
number theory (cf. \cite{BI}, \cite{S}). In spite of the fact that
the whole theory was originally motivated by problems in complex
analysis (H.Bohr \cite{Bo}, p.3), analytic aspects of the theory
are less known. Holomorphic almost periodic functions have certain
specific properties, mainly because almost periodicity of a
holomorphic function causes strong restrictions on the
distribution of its values (for example, of its zeroes). The main
contributions to the classical theory of holomorphic almost
periodic functions of one variable are due to B.\,Jessen,
H.Tornehave, B.\,Ja.\,Levin and M.\,G.\,Krejn; for a detailed
presentation of the subject, see \cite{JT} and \cite{L}.

The case of several complex variables is much more difficult.
First results in this direction concerned zero sets of exponential
sums (O.\,Gelfond \cite{G1}, \cite{G2},  B.\,Ja.\,Kazarnovskii
\cite{K1}, \cite{K2}). A fruitful approach developed by
L.\,I.\,Ronkin has allowed  to investigate general holomorphic
almost periodic functions and mappings (cf.\cite{Ra1}--\cite{R3}).
This resulted, in particular, in the notions of almost periodic
divisors and holomorphic chains (cf. \cite{R5}, \cite{R6},
\cite{FRR2}, \cite{FRR3}; also see the survey \cite{FR}). Note
that the main accent in those papers was made to the asymptotic
behavior and asymptotic characteristics.

The problem of inner description for the divisors of holomorphic
almost periodic functions has been raised by M.\,G.\,Krejn and B.\,Ja.\,Levin
in \cite{KL}, where it was solved for entire  almost periodic functions of
exponential type on the plane with zeroes in a strip of finite width.
Note that the divisors of holomorphic almost periodic functions
inherit, in a sense, the property of being almost periodic,
nevertheless there exist almost periodic divisors that are not the divisors of any almost
periodic holomorphic functions (cf. \cite{T2} and \cite{R6}).
For holomorphic almost periodic functions on a strip with spectrum in a free group
the problem was solved in \cite{T1}. A complete description of the divisors of
holomorphic and meromorphic almost functions on a strip was
given in \cite{F} and \cite{FP}.

Notice that almost periodic divisors in a strip on the plane were
classified in \cite{T1} with the help of certain  integer valued
matrix, while in \cite{F}, \cite{FP} this was done in terms of
cohomology classes from $H^2(K,\Z)$, where $K$ is the so-called
Bohr's compactification of the real axis (see Section 1).

The multi-dimensional problem of realizability of almost periodic
divisors as the divisors of almost periodic holomorphic functions
was studied by L.\,I.\,Ronkin  for periodic divisors
(cf.\cite{R4}) and for their restrictions to certain family of
planes (cf.\cite{R6}) with the help of an integer valued matrix,
too. But for arbitrary almost periodic divisors the problem has
never been raised before.

Our approach to the realizability problem is very close to the
classical method of investigation of the Second Cousin Problem on
a domain $D\subset\C^m$ (see \cite{H1}, Ch.5). Namely, in the classical case
a line bundle over $D$ corresponds to each data for the Cousin Problem on $D$
such that the problem has a solution if and only if this
bundle is trivial. Therefore, the first Chern class (i.e., the
corresponding element of the group $H^2(D,\Z)$) of the bundle
is assigned to each data such that the problem is solvable if and
only if this class is trivial. In particular, if $H^2(D,\Z)=\{0\}$,
then every Second Cousin Problem on $D$ has a solution.
In our investigation, a line bundle over Bohr's compact set $K$ is assigned
to each almost periodic divisor on
a tube domain $\T_\O$ (the case $T_\O=\C^m$ is not excluded),
such that the bundle is trivial just for the divisors
of holomorphic almost periodic functions on the domain;
therefore, the first Chern class of the bundle corresponds to each
almost periodic divisor in the domain such that it is the divisor of a holomorphic
almost periodic function if and only if this class is trivial.
Note that $H^2(K,\Z)\neq\{0\}$, therefore some problems of realizability
for almost periodic divisors have no solutions.

As for the case of the Second Cousin Problem, we need to solve the
appropriate $\bar\partial$-problem; we use a technique from \cite{BA} to get
a required integral representation for a solution with specific properties of
the $\bar\partial$-problem in a tube domain.

The paper is organized as follows.

In Section 1, we give the main definitions, necessary
notations, and some information about almost periodic
functions and bundles.

In Section 2, we introduce the notion of holomorphic function on
$K\times\O$ and establish a correspondence between such functions
and holomorphic almost periodic functions in $T_\O$. Then we
construct a line bundle over $K\times\O$, corresponding to an
almost periodic divisor in $T_\O$.

In Section 3, we solve the
appropriate $\bar\partial$-problem on $\T_\O$ and prove that if the
constructed bundle is trivial, then the divisor is the divisor of a
holomorphic almost periodic function in $T_\O$.

In Section 4, we show that the Chern class of the bundle is
actually an element of the group $H^2(K,\Z)$ and obtain simple
geometric conditions sufficient for realizability of divisors as
the divisors of holomorphic almost periodic functions. For the
case of divisors with spectrum in a finitely-generated additive
subgroup of $\R^m$, we establish a correspondence between Chern
classes of the divisors and skew-symmetric matrices with integer
entries.

In Section 5, we investigate the Chern classes of periodic
divisors and some classes of almost periodic divisors.
We also find a structure formula for Chern classes of almost
periodic divisors and prove that an arbitrary Chern class
is a finite sum of the Chern classes for periodic divisors.

In Section 6, we obtain a complete classification of the divisors
of homogeneous coordinates of almost periodic holomorphic mappings from tube
domains into projective spaces.

\bigskip

\section{Definitions, notations, and some preliminary information}
\label{sec1}

\bigskip

A continuous function $f(x)$ on $\R^m$ is called almost periodic
if the collection $\{S_tf\}_{t\in\R^m}$ is a relatively compact
set with respect to the topology of uniform convergence on $\R^m$; here
$S_tf(x)=f(x+t)$ is the shift along $t$.

The class of such functions coincides with the closure, with
respect to the topology of the uniform convergence on $\R^m$, of
the set of all finite exponential sums \be \label{sum1} \sum
a_ne^{i\<x,\l_n\>},\quad \l_n\in\R^m, \en where
$\<x,\l_n\>=x^1\l_n^1+\dots+x^m\l_n^m, \quad a_n\in\C$. In the
case $m=1$, this is one of the main results of the classical
theory of almost periodic functions.

A continuous function $f(z)$ on a tube domain
$T_\O=\{z=x+iy:\,x\in\R^m,\,y\in\O\subset\R^m\}$ is called almost
periodic on $T_\O$ if the collection $\{S_tf\}_{t\in\R^m}$ is a relatively
compact set with respect to the topology of uniform convergence on every
tube domain $T_{\O'}$, $\O'\Subset\O$.

The class of such functions coincides with the closure of the set of
all finite exponential sums (\ref{sum1}), where $a_n$ are
continuous functions on $\O$.

Throughout the paper, $AP(\R^m)$ is the class of all almost periodic
functions on $\R^m$, $AP(T_\O)$ is the class of all almost periodic
functions on $T_\O$, $H(D)$ is the class of all holomorphic
functions on $D$, $APH(T_\O)=AP(T_\O)\cap H(T_\O)$.
We always consider tube domains with convex base $\Omega$.

Note that every function $f\in APH(T_\O)$ belongs to the closure
of the set of all finite sums
\be
\label{sum2} \sum a_ne^{i\<z,\l_n\>},\quad \l_n\in\R^m,\quad a_n\in\C.
\en

The notion of spectrum for functions $f$ from $AP(\R^m)$ or
$AP(T_\O)$ is introduced as
\be \label{spec} \sp f
=\{\l\in\R^m:\,\lim_{T\to\infty}(2T)^{-m}
\int_{[-T,T]^m}f(x+iy)e^{-i\<x,\l\>}dx\not\equiv 0\}; \en the
spectrum is at most countable, and all the exponents $\l$ in sums
(\ref{sum1}) or (\ref{sum2}) belong to the  spectrum of $f$.

For the case of functions on the axis or on a strip, all these
facts are basic results of the theory of almost periodic function
(see, for example, \cite{C}); for the general case the proves are
similar, see \cite{R3}, \cite{R5}.

By $C^\infty_{p,q}(T_\O)$ we denote the space of $(p,q)$-forms
\be\label{sum3} \sum
a_{i_1,\dots,i_p,j_1,\dots,j_q}(z)dz^{i_1}\wedge\dots\wedge
dz^{i_p} \wedge d\bar z^{j_1}\wedge\dots\wedge d\bar z^{j_q} \en
with infinitely smooth coefficients
$a_{i_1,\dots,i_p,j_1,\dots,j_q}(z)$; we assume that the
coefficients and all their derivatives are bounded on any domain
$T_{\O'}$ for $\O'\Subset\O$. Also, by $C^0_{p,q}(T_\O)$ we denote
the space of $(p,q)$-forms such that all the coefficients are
continuous and have compact supports. Finally, by $M_{p,q}(T_\O )$
we denote the space of $(p,q)$-currents of order $0$ on $T_\O$,
i.e., forms (\ref{sum3}) whose coefficients are complex measures
on $T_\O$. This space is equipped with the weak topology of
convergence on elements of $C^0_{m-p,m-q}(T_\O)$.

A typical example of a current from $M_{1,1}(T_\O)$ is the
current of integration over the divisor $d_F$ of a function $F\in
H(T_\O)$, i.e., the current $\dd\log|F(z)|$. We will identify
occasionally a divisor $d$ with the corresponding current of
integration and write $d\in M_{1,1}(T_\O)$.

A current $\theta\in M_{p,q}(T_\O)$ is almost periodic if for any
$\phi\in C^0_{m-p,m-q}(T_\O)$ the function $(\theta,S_t\phi)$, as
a function of $t$, belongs to $AP(\R^m)$; the space of such
currents is denoted by $APM_{p,q}(T_\O)$. A divisor $d$ is almost
periodic if $d\in APM_{1,1}(T_\O)$; all the divisors of functions
$\fHT$ are almost periodic (cf. \cite{R5}); if the measure
$\sum_k\partial^2\log|F|/\partial z_k \bar\partial z_k$ for
$F\in\T_\O$ is almost periodic, then the divisor $d_F$ is almost
periodic (see \cite{FRR3}).

We will make use of the following result (cf.\cite{R5},
\cite{FR}).

\begin{prop}\label{t-con}
If a sequence $\{f_n\}\subset H(D)$ converges uniformly
on every compact subset of $D$, then the sequence $\{\log |f_n|\}$
converges on $D$ in the sense of distributions.
\end{prop}

It was proved in \cite{R3} that the coefficients of a current
$\theta\in APM_{p,q}(T_\O)$ can be approximated in the weak
topology on sets $\{S_t\phi\}_{t\in\R^m}$, $\phi\in
C^0_{0,0}(T_\O)$, uniformly in $t\in\R^m$, by finite sums
(\ref{sum1}); in this case $a_n$ are complex measures on $\O$.

The spectrum of $\theta\in APM_{p,q}(T_\O)$ is the union of the
spectra of the functions $(\theta,S_t\phi)$ over all $\phi\in
C^0_{m-p,m-q}(T_\O)$. This set is countable as well, and all the
exponents in sum (\ref{sum1}) belong to the  spectrum of $\theta$
(see \cite{R5} or \cite{FR}).

For any additive subgroup $\G\subset\R^m$, we denote by $AP(\R^m,\G)$,
$AP(T_\O,\G)$, $APH(T_\O,\G)$, and $APM_{p,q}(T_\O,\G)$ the
corresponding classes of functions with spectra in $\G$.

Bohr's compactification $\hR$ is the closure of the image of
$\R^m$ under the map $j$ into Tikhonov's product of the circles
$\T_\l$, \be \label{j}
j:\,x\mapsto\{e^{i\<x,\l\>}\}\in\prod_{\l\in\R^m}\T_\l. \en The
preimage of Tikhonov's topology on $\R^m$ with respect to the map
$j$ is called Bohr's topology and is denoted by $\beta$. It is
clear that this topology is the weakest one where all the sums
(\ref{sum1}) are continuous. If we replace  the product of the
circles over all $\l\in\R^m$ by the product over $\l\in\G$, $\G$
being a subgroup of $\R^m$, then the corresponding Bohr's
compactification is denoted by $\hG$ and the corresponding
topology by $\beta(\G)$. Evidently, one can take the product of
the circles only over all generators of $\G$. In particular, for
$\G=\Lin_\Z\{\l_1,\dots,\l_N\}$ with vectors
$\l_1,\dots,\l_N\in\R^m$ linearly independent over $\Z$, the set
$\hG$ coincides with the torus $\T^N$. Then it follows from the
approximation property that functions from $AP(\R^m)$ are
uniformly continuous with respect to topology $\beta$ and
functions from $AP(\R^m,\G)$ are uniformly continuous with respect
to topology $\beta(\G)$. Hence functions from $AP(\R^m)$ are just
the restrictions to $j\R^m$ of continuous functions on $\hR$, and
functions from $AP(\R^m,\G)$ are just the restrictions to $j\R^m$
of continuous functions on $\hG$ (see \cite{AK} or \cite{B}).

Following \cite{H1} and \cite{H}, we will consider line bundles $\F$
with the fiber $\C$ over a
paracompact base $X$; in this paper $X=\hR$ or $X=\hRO$, where
$\O$ is an open convex subset of $\R^m$. We say that a section $\p$
of the bundle over an open set $\o\subset\hRO$ is holomorphic if
$\p(\tau\cdot jx,y)$ is a smooth function of $x,\,y\in\R^m$ for
small $x$, and $\bar\partial_z \p(\tau\cdot jx,y)\mid_{x=0}=0$
for all $(\tau,y)\in\o$, $j$ being defined in (\ref{j}). To each
line bundle $\F$ over $X$ assign the (first) Chern  class
$c(\F)\in H^2(X,\Z)$. Namely, let $\o_\a$ be a sufficiently fine
covering of $X$ and $g_{\a,\beta}(x),\ x\in X$, be the
corresponding transition functions. The integer
\be \label{coh}
c_{\a,\beta,\g}= {1\over 2\pi i}\log g_{\a,\beta}(x)+{1\over 2\pi
i}\log g_{\beta,\gamma}(x) +{1\over 2\pi i}\log g_{\gamma,\a}(x)
\en
corresponds to every triple of indices such that
$\o_\a\cap\o_\beta\cap\o_\gamma\neq\emptyset$ (if
$\o_\a\cap\o_\beta\cap\o_\gamma=\emptyset$, take
$c_{\a,\beta,\gamma}=0$). The collection
$\{c_{\a,\beta,\gamma}\}$ forms a 2-cocycle; its cohomology class
$[\{c_{\a,\beta,\gamma}\}]\in H^2(\{\o_\a\},\Z)$ does not depend
on the choice of branch of logarithm in (\ref{coh}). Moreover, if
we multiply each function $g_{\a,\beta}(x)$ by
$h_\a(x)/h_\beta(x)$, $h_\a(x)\neq 0$ for all $\a$ and
$x\in\o_\a$, then the cohomology class remains the same. Passing
to the inductive limit with respect to refinement of covering, we obtain
the first Chern class $c(\F)\in H^2(X,\Z)$. If there exists a
global section $\{\p_\a\}$ vanishing nowhere, then we can take
$\log g_{\a,\beta}=\log\p_\a-\log\p_\beta$. This yields that
$c_{\a,\beta,\g}=0$, therefore $c(\F)=0$.

On the other hand, to each $c\in H^2(X,\Z)$ one can assign an
element of the group $H^1(X,\Xi^*)$, where $\Xi^*$ is the
multiplicative bundle of germs of continuous functions vanishing
nowhere on $X$. Namely, the 1-cocycle
$\{g_{\a,\beta}\}$ of the bundle $\Xi^*$ corresponds to
$c_{\a,\beta,\g}\in H^2(\{\o_\a\},\Z)$ such that (\ref{coh}) is
fulfilled. Note that the 1-cocycles $\{g_{\a,\beta}\}$ satisfy the
conditions
$$
g_{\a,\beta} g_{\beta,\a}=1\quad{\rm and}\quad
g_{\a,\beta}g_{\beta,\g}g_{\g,\a}=1,
$$
i.e., $\{g_{\a,\beta}\}$ are transition functions for some line
bundle over $X$.

It is clear that the line bundle with the transition functions
$g_{\a,\beta}/\tilde g_{\a,\beta}$ corresponds to the difference
of cohomology classes $c-\tilde c$. Further, if $c=0$, then
$c_{\a,\beta,\g}=0$ for sufficiently fine covering. Now it follows
in the standard way that $g_{\a,\beta}=\p_\a/\p_\beta$ for some
global section $\{\p_\a\}$ vanishing nowhere. Note also that if
the Chern class $c$ corresponds to the bundles $\F$ and $\tilde\F$
with transition functions $\{g_{\a,\beta}\}$ and $\{\tilde
g_{\a,\beta}\}$, respectively, then the bundle with the transition
functions $\{g_{\a,\beta}/\tilde g_{\a,\beta}\}$ has a global
section vanishing nowhere.

\bigskip

\section{Almost periodic divisors and line bundles}
\label{sec2}

\bigskip

First we give some auxiliary results.

\begin{prop}
\label{t-ex1} Each function $\fAT$ extends to the set $\hRO$ as a
continuous function $\f$ such that $\bar f(jx,y)=f(z)$.
Conversely, each continuous function $\f$ on $\hRO$ generates the
function $f(z)=\bar f(jx,y)\in AP(T_\O)$. In addition, $\fHT$ iff
\be \label{hol} \bar\partial_z\fxy\mid_{x=0}=0 \en for all
$(\tau,y)\in\hRO$.
\end{prop}

\ws It is obvious that every function $e^{i\<x,\l\>}$ extends to
$\hR$ as a continuous function $\bar e_\l(\tau)$, therefore every
exponential sum (\ref{sum1}) with coefficients $a_n(y)$ continuous
on $\O$ extends to $\hRO$ as a continuous function, too. The same
assertion is true for every $\fAT$ as a uniform limit of such sums
on every $T_{\O'},\,\O'\Subset\O$. Since $f(z)=\bar f(jx,y)$ is
valid for the functions $e^{i\<x,\l\>}$, the equality holds for
all $\fAT$. On the other hand, if $\bar f(\tau,y)$ is continuous
on $\hRO$, then $\bar f(\tau,y)$ is uniformly continuous on
$\hR\times P$ for each compact set $P\subset\O$, therefore the
family of shifts $\{S_tf\}_{t\in\R^m}$ for functions $f(z)=\bar
f(jx,y)$ is relatively compact set with respect to the topology of
uniform convergence on every tube domain $T_{\O'}$,
$\O'\Subset\O$, and $\fAT$. Finally, since $j\R^m$ is dense in
$\hR$ and $\bar f(j(t+x),y)=S_tf(z)$, the function $\fxy$ is
holomorphic, too.\bs

\begin{prop}
\label{t-ex2} Each divisor $d\in APM_{1,1}(T_\O)$ defines a
unique continuous mapping $\bar d:\,\hR\mapsto M_{1,1}(T_\O)$
with the properties
\be \label{div1}
\bar d(\tau\cdot jt)=S_t\bar d(\tau),
\en
\be \label{div2}
\bar d(1)=d.
\en
Conversely, each
continuous mapping $\bar d$ from $\hR$ to $M_{1,1}(T_\O)$ with
the property (\ref{div1}) generates a divisor $d\in
APM_{1,1}(T_\O)$ by equality (\ref{div2}); if $d$ is the divisor
of a function $\fHT$, then $\bar d(\tau)$ is the divisor of the
holomorphic function $\bar f(\tau\cdot jx,y)$.
\end{prop}

\ws Let $\bar e_\l(\tau)$ be the function from the proof of the
previous proposition. Then for every $\l\in\R^m$ the map
$\tau\mapsto\bar e_\l(\tau\cdot jx)$ is a
continuous map from $\hR$ to $AP(\R^m)$  with the property
(\ref{div1}). Therefore every sum (\ref{sum1})
with complex measures $a_n(y)$ on $\O$ generates a continuous map
from $\hR$ to $M_{1,1}(T_\O)$. Every divisor $d\in APM_{1,1}$ can
be approximated by such sums, hence it generates a continuous map
$\bar d:\,\hR\mapsto AP(T_\O)$  with the property (\ref{div1}).
Conversely, each continuous map $\bar d$ from
$\hR$ to $AP(\R^m)$  with the properties (\ref{div1}) and each
form $\a\in C^\infty_{1,1}(T_\O)$ induce the continuous function
$(\bar d(\tau),\a)$ on $\hR$. Hence the function $(\bar
d(1),S_t\a)=(\bar d(jt),\a)$ belongs to $AP(R^m)$. This means that
$d\in APM_{1,1}(T_\O)$. Using  the density of $jR^m$ in $\hR$ and
the equality $\bar d(jt)=S_td$, we get the uniqueness of $\bar d$.
Finally, if $d=d_f$ for $\fHT$, then the map
$\hat d:\,\tau\mapsto\dd\log|\fxy|$ satisfies (\ref{div1}) and
(\ref{div2}); this map is continuous because the map
$\tau\mapsto\fxy$ is continuous as a map from $\hR$ to
$APH(T_\O)$ and, by Proposition \ref{t-con}, the map $f\mapsto\log|f|$
is continuous as a map from $H(T_\O)$ to $M_{0,0}(T_\O)$.
Thus $\bar d=\hat d$.\bs
\medskip

{\bf Remark.} In the previous propositions, the classes
$AP(T_\O)$, $APM_{1,1}(T_\O)$ and the compact set $\hR$ can be
replaced by $AP(T_\O,\G)$, $APM_{1,1}(T_\O,\G)$ and the compact
set $\hG$, respectively. Really, we can use for approximation of
the functions and currents with spectrum in $\G$ only the sums of
exponents $e^{i\<x,\l\>}, \ \l\in\G$. Conversely, the restriction
to the set $j\R^m\times\O$ of a continuous on $\hG\times\O$
function $\bar f(\tau,y)$ is  a continuous function $f(z)$ on $x$
with respect to the topology $\beta(\G)$, therefore $\sp
f(x+iy_0)\subset\G$ for each fixed $y_0\in\O$, and we obtain $\sp
f\subset\G$. In the same way,  we get $\sp(\bar d(jt))\subset\G$
in Proposition~\ref{t-ex2}.

\medskip

We say that a function $\p(\tau,y)$, continuous on an open set
$\o\subset\hR\times\O$,
 is holomorphic on $\o$ if the function
$\pxy$ is smooth for small $x$ and satisfies equation (\ref{hol})
for all $(\tau,y)\in\o$. Further, $\p(\tau,y)$ defines the divisor
$d\in APM_{1,1}(T_\O)$ on $\o$ if for each $\tau,y_0\in\o$ there
exists $\d>0$ such that the restriction of the divisor $\bar
d(\tau)$ to the ball $B(iy_0,\d)\subset\C^m$ coincides with the
restriction of the divisor of the function $\pxy$.

It can be shown that a holomorphic function defining a divisor is
essentially unique. Namely, we have

\begin{prop} \label{t-quot}
Suppose $\bar f(\tau,y),\ \bar g(\tau,y)$ are
holomorphic functions on $\o\subset\hRO$. If for each point
$z_0=x_0+iy_0,\ (jx_0,y_0)\in\o$, the quotient $\bar f(jx,y)/\bar
g(jx,y)$ is a holomorphic function on some neighborhood of $z_0$,
then $\bar f(\tau,y)=\bar g(\tau,y)\bar h(\tau,y)$ with
holomorphic on $\o$ function $\bar h(\tau,y)$. Moreover, if the
functions $\bar f(jx,y)$ and $\bar g(jx,y)$ have the same divisor,
then $\bar h(\tau,y)\neq 0$.
\end{prop}

\ws Fix $(\tau_0,y_0)\in\o$. If $\bar g(\tau_0,y_0)\neq 0$, then we take
$\bar h=\bar f/\bar g$ on some neighborhood of $(\tau_0,y_0)$. Suppose
$\bar g(\tau_0,y_0)=0$. Let $z=x+iy=(x^1,'x)+i(y^1,'y),\,
x^1,y^1\in\R,\,'x,'y\in\R^{m-1}$.
Without loss of generality it can be assumed that
$a(z^1)=\bar g(\tau_0\cdot j(x^1,'0),(y_0^1,'y_0))\not\equiv 0$. We have
$|a(z^1)|\ge\gamma>0$ whenever $|z^1-iy_0^1|=\e$ and $\e$ is sufficiently
small.  Let $\tau$ be near $\tau_0$, $'y$ be near $'y_0$,
$\zeta=\xi+i\eta\in\C$. Take
$$
\bar h(\tau,y)={1\over 2\pi i}\int_{|\zeta-iy^1_0|=\e}
{\bar f(\tau\cdot j(\xi,'0),(\eta,'y))\over
\bar g(\tau\cdot j(\xi,'0),(\eta,'y))}{d\zeta\over \zeta-iy^1}.
$$
This function is continuous in some neighborhood of
$(\tau_0,y_0)$. Take $\tau=jt,\ t\in\R^m$. Using Cauchy's integral
formula, we obtain $\bar h(jt,y)\bar g(jt,y)=\bar f(jt,y)$. Now
the first statement of the proposition follows from the density
of $\R^m$ in $\hR$. To prove the last statement we interchange
$\bar f$ and $\bar g$. \bs

\begin{prop}
\label{t-loc}
For any divisor $d\in APM_{1,1}(T_\O)$ and any point
$(\tau_0,y_0)\in\hRO$ there exists a neighborhood $\o$ and
a holomorphic function $r(\tau,y)$ defining the divisor $d$ on $\o$.

\end{prop}

\ws Let $\bar d(\tau)$ be  the continuous map from proposition
\ref{t-ex2}. Since $T_\O$ is a convex domain, we see that for each
$\tau\in\hR$ there exists a function $f_\tau(z)\in H(T_\O)$ such
that $\dd\log|f_\tau(z)|$ is the current of integration over the
divisor $\bar d(\tau)$. Let
$z=(z^1,'z)=(x^1,'x)+i(y^1,'y),\,x^1,y^1\in\R,\,'x,'y\in\R^{m-1},\,
\zeta=z^1=x^1+iy^1$. Using a linear transformation of coordinates
with real coefficients, we may assume that
$f_{\tau_0}(\zeta,i'y_0)\not\equiv 0$. Fix $\d>0$ such that
$f_{\tau_0}(\zeta,i'y_0)\neq 0$ for $0<|\zeta-iy_0^1|\le\d$.
Therefore we have $\supp\bar d(\tau_0)\cap\{z=(\zeta,i'y_0)\in
T_\O:\, \d/2\le|\zeta-iy^1_0|\le\d\}=\emptyset$. Hence there
exists a neighborhood $U_{\tau_0}$ of $\tau_0$ and $\e>0$ such that for
$\tau\in U_{\tau_0},\ |'z-i'y_0|<\e$ the functions
$f_\tau(\zeta,'z)$ as functions of $\zeta$ have no zeros on the
set $\d/2\le|\zeta-iy^1_0|\le\d$.

Let $\zeta_1(\tau,'z),\dots,\zeta_k(\tau,'z)$, $k=k(\tau,'z)$, be
all zeros of the function $f_\tau(\zeta,'z)$ on the set
$|\zeta-iy^1_0|<\d/2$ as a function of $\zeta$. Let $\p(\zeta)$
be a nonnegative infinitely smooth function such that
$\p(\zeta)=1$ for $|\zeta-iy_0^1|\le\d/2$ and $\p(\zeta)=0$ for
$|\zeta-iy_0^1|\ge\d$, $\Delta_\zeta$ be Laplace operator on the
plane $C_\zeta$. We have
$$
k(\tau,'z)=2\pi(\Delta_\zeta\log|f_\tau(\zeta,'z)|,\p(\zeta))=\pi/i
\int\log|f_\tau(\zeta,'z)|\Delta_\zeta\p(\zeta)d\zeta\wedge d\bar \zeta.
$$
It follows from Hurwitz' Theorem that for any fixed $\tau\in
U_{\tau_0}$ the number $k(\tau,'z)$ does not depend on $'z$ from
the disc $|'z-i'y_0|<\e$.

Let $\psi('z)$ be an infinitely smooth function with support in the set
$\{'z:\,|'z-i'y_0|<\e\}$ such that $(2i)^{1-m}\int \psi('z)
d'z\wedge d\bar{'z}=1$. Then
$$
k(\tau,'z)=2\pi(2i)^{-m}\int\log|f_\tau(\zeta,'z)|\psi('z)
\Delta_\zeta\p(\zeta)d\zeta\wedge d\bar \zeta\wedge d'z\wedge d\bar{'z}.
$$
The integral is the action of the current
$\dd\log|f_\tau|$ on the form $\psi('z)\p(\zeta)d'z\wedge
d\bar{'z}\in c^\infty_{m-1,m-1}(T_\O)$, hence it is continuous as
a function of $\tau$, and $k(\tau,'z)\equiv\const$ on
$U_{\tau_0}\times\{'z:\ |'z-i'y_0|<\e\}$. Since $\d$ is arbitrary
small, the functions $\zeta_1(\tau,'z),\dots,\zeta_k(\tau,'z)$
are continuous at the point $(\tau_0,i'y_0)$. By the same
arguments, these functions are continuous at each point
$(\tau,'z)$ with the property $f_\tau(\zeta,'z)\not\equiv 0$.

Now we define a function $r(\tau,y)$. If $iy_0\not\in\supp\bar d(\tau_0)$,
we can take $r(\tau,y)\equiv 1$ and $\o$ is small enough. Otherwise we
take
$\o=\{(\tau,y):\,\tau\in U_{\tau_0},|y^1-y^1_0|<\d,|'y-'y_0|<\e\},\
r(\tau,y)=P(\tau,iy^1,i'y)$
with
$P(\tau,\zeta,'z)=(\zeta-\zeta_1(\tau,'z))\cdots(\zeta-\zeta_k(\tau,'z))$.
Represent this polynomial in the form
$P(\tau,\zeta,'z)=\zeta^k+b_{k-1}(\tau,'z)\zeta^{k-1}+\dots+b_0(\tau,'z)$.
Since
$$
\sum_{1\le j\le k}(\zeta_j(\tau,'z))^s={1\over 2\pi i}
\int_{|\zeta-iy^1_0|=\d}{\partial f_\tau(\zeta,'z)\over \partial\zeta}
{\zeta^sd\zeta\over f_\tau(\zeta,'z)},
$$
all the functions $b_l(\tau,'z),\ l=1,\dots,k-1,$ are holomorphic
on the set $\{'z:\,|'z-i'y_0|<\e\}$. Since the polynomial
$P(\tau,z^1,'z)$ has the same zeros as $f_\tau(z)$ on the set
$\{z:\,|z^1-iy^1_0|<\d,\,|'z-'y_0|<\e\}$, we get that the
restriction of the divisor of $P(\tau,z)$ to this set is $\bar d(\tau)$.
It follows from
(\ref{div1}) that $P(\tau\cdot jt,z)=P(\tau,t+z)$ for small $t$ and $z$,
hence $r(\tau,y)$ is a holomorphic function on $\o$. \bs

\medskip

Now we introduce the main result of this section.

\begin{th}
\label{t-bund} To each $d\in\AM$ one can assign a line bundle
$\F_d$ over $\hRO$ and a global holomorphic section of this bundle
defining the divisor $d$; $d$ is the divisor of a function $\fHT$
iff $\F_d$ has a global holomorphic section vanishing nowhere.
\end{th}

\ws Take a sufficiently fine covering $\{\o_\a \}$ of $\hRO$
such that the divisor $d$ is defined by a holomorphic function
$r_\a(\tau,y)$ on each $\o_\a$. Taking into account Proposition \ref{t-quot},
we get that the functions $g_{\a,\beta}(\tau,y)=r_\a(\tau,y)/r_\beta(\tau,y)$
vanish nowhere on $\o_\a\cap\o_\beta$. It is obvious that these functions
satisfy all the conditions for being transition functions of some line bundle
$\F_d$ over $\hRO$. The equality
\be
\label{trans}
r_\beta(\tau,y)=r_\a(\tau,y) g_{\a,\beta}(\tau,y)
\en
means that
$\{r_a\}$ is the section of $\F_d$. If $d=d_f$ for some $\fHT$,
then the function $\bar f(\tau,y)$ from Proposition \ref{t-ex1} also
defines $d$,
hence $\{r_a(\tau,y)/\bar f(\tau,y)\}$ is a global holomorphic section
vanishing nowhere.

Conversely, if $\{R_\a(\tau,y)\}$ is a global holomorphic section of $\F_d$
vanishing nowhere, then the holomorphic function $\bar F(\tau,y)=
r_a(\tau,y)/R_a(\tau,y)$ is well defined on $\hRO$ and defines the same
divisor as $r_a(\tau,y)$. \bs

\medskip

Using Remark to Propositions \ref{t-ex1} and \ref{t-ex2} and
arguing as above, we also get an analog of this theorem for
divisors with spectrum in a subgroup $\G$:

\begin{th}
\label{t-bundG} To each $d\in APM_{1,1}(T_\O,\G)$ one can assign a
line bundle over $\hG\times\O$ such that some global holomorphic
section of this bundle defines the divisor $d$; $d$ is the
divisor of a function $f\in APH(T_\O,\G)$ iff the bundle has a
global holomorphic section vanishing nowhere.
\end{th}

In the end of this section we prove the converse statement to
Theorem \ref{t-bund}:

\begin{prop}
\label{t-bund1}
Every holomorphic section $\{r_a(\tau,y)\}$ of a line bundle $\F$
over $\hRO$ defines some divisor $d\in\AM$.
\end{prop}

\ws Fix $\tau\in\hR$. The equality
\be\label{cur}
\bar d(\tau)=\dd\log|r_\a(\tau\cdot jx,y)|
\en
for $\a$ such that $(\tau\cdot jx_0,y_0)\in\o_\a$
defines the divisor $\bar d(\tau)$ on the set
$\{z=x+iy:\,|z-z_0|<\e\}$, $\e$ is small enough. It follows from
(\ref{trans}) that the current (\ref{cur}) is well-defined and belongs to
$M_{1,1}(T_\O)$. Since $(\tau\cdot jt)\cdot jx_0=\tau\cdot
j(x_0+t)$, we obtain (\ref{div1}). Furthermore, the function
$r_\a(\tau\cdot jx,y)$ is continuous with respect to $\tau$; using
Proposition \ref{t-con}, we obtain that the function $\log|r_\a(\tau\cdot
jx,y)|$, as a mapping from $\hR$ to $M_{0,0}(\{z:\,|z-z_0|<\e\}$,
is continuous. Therefore the mapping $\bar d(\tau)$ of $\hR$ to
$M_{0,0}(T_\O)$ is continuous, and we can use Proposition
\ref{t-ex2}. \bs
\bigskip

\section{Existence of a global holomorphic section vanishing nowhere}
\label{sec3}

\bigskip
Here we get a sufficient condition for the existence of a global
holomorphic section vanishing nowhere. First we need the following
result.

\begin{prop}
\label{t-d} For each pair of convex domains $\O'\Subset\O$ the
equation $\bp u=a$ with a $\bp$-closed $(0,1)$-differential form
$a(z)=\sum_k a_k(z)\bar dz^k$ with bounded and $C^1$-smooth on
$T_\O$ coefficients has a solution $Aa\in C_{0,0}^\infty(T_{\O'})$
with the following properties:

a) $A$ is a linear operator,

b) $AS_ta=S_tAa$ for all $t\in\R^m$,

c) $\sup_{z\in T_{\O'}}|Aa(z)|\le c(\O,\O')\sup_{z\in T_\O}\left(\sum_k|a_k(z)|^2\right)^{1/2}$.
\end{prop}

\medskip

\ws At first, assume that $\O$ is bounded and
$\O=\{y:\,\kappa(y)<0\}$ for some convex and $C^2$-smooth function
$\kappa(y)$ in a neighborhood of $\O$,
$\grad_y\kappa(y)\neq 0$ on $\partial\O$, and $a(z)$ is any $\bp$-closed
$(0,1)$-differential form with bounded and $C^1$-smooth on $T_{\overline{\O}}$
coefficients. Let $\beta(t)$
be a convex, $C^2$-smooth function on $\R$ such that $\beta(t)=0$
for $t<0$, $\beta(t)>0$ for $t>0$, and $\beta(1)>-\min_{\O}\kappa(y)$.
For $z=x+iy\in\C^m$ and $r>0$, put
$\rho(z)=\kappa(y),\; \rho_r(z)=\rho(z)+ \beta(|x|-r),\ G_r=\{z:\,\rho_r(z)<0\}$.
It is clear that $\{z\in
T_\O:\,|x|<r\}\subset G_r\subset\{z\in T_\O:\,|x|<r+1\}$.

Let $\zeta=\xi+i\eta\in\C^m,\
\<z,\zeta\>=\sum_{k=1}^m z^k\zeta^k$. Fix a compact set
$K\subset T_\O$. We claim that the inequality
\be
\label{grad1} |\<\grad_\zeta\rho_r(\zeta),\zeta-z\>|\ge\e
\en
is true with some $\e>0$ for all $z\in K,\ \zeta\in\T_{\bar\O}$ such that
$\dist(\zeta,\partial G_r)<\d,\ \d=\d(K)$, for all $r>r(K)$.
It is sufficient to estimate from below the sum \be
\label{grad2}
\<\grad_\xi\beta(|\xi|-r),\xi-x\>+\<\grad_\eta\kappa(\eta),\eta-y\>. \en
Here the first term is nonnegative for $|x|<r$ and all
$\xi\in\R^m$. Since
\be\label{grad3}
\<\grad_\eta\kappa(\eta),\eta-y\>\ge 2\e>0
\en
for all
$\eta\in\partial\O,\ z\in K$, and some $\e>0$, we see that the
value (\ref{grad2}) is at least $\e$ whenever
$\dist(\eta,\partial\O)\le 2\d$ for some $\d>0$. On the other
hand, if $\dist(\zeta,\partial G_r)<\d$ and $\dist(\eta,\partial\O)>2\d$, then
$|\xi|-r\ge\gamma$ for some $\gamma>0$, and  we have
$$
\<\grad_\xi\beta(|\xi|-r),\xi-x\>=\beta'(|\xi|-r)\<\xi,\xi-x\>/|\xi|\ge
\beta'(\gamma)(|\xi|-|x|).
$$
Since the last term is greater than $\beta'(\gamma)(r+\gamma-|x|)$, and
since $|\<\grad_\eta\kappa(\eta),\eta-y\>|\le c<\infty$ for all
$\eta,y\in\O$, we get the desired estimate for a sufficiently large $r$.

Furthermore, put \be\label{BA1}
u_r(z)=C_m\int_{G_r}\exp\<z-\zeta,\,\zeta-z\>
{a\wedge\<s_r,\,d\zeta\>\wedge(d_\zeta\<s_r,\,d\zeta\>)^{m-1}\over
\<s_r,\,\zeta-z\>^m} \en with $s_r=\overline{\<\grad_\zeta
\rho_r(\zeta),\,\zeta-z\>}
\grad_\zeta\rho_r(\zeta)+|\rho_r(\zeta)|(\bar\zeta-\bar z)$. Note
that the kernel in (\ref{BA1}) is, up to a constant, just the
component of bidegree $(0,\,0)$ in $z$ of the kernel given by
equality (14) of \cite{BA} for $Q=z-\zeta$. Moreover, it follows
from (\ref{grad2}) that $s_r$ satisfies the necessary conditions
of Theorem 2 of \cite{BA} (see also item ii) on p.~102 there),
hence for the appropriate choice of the constant $C_m$ in
(\ref{BA1}) we have $\bp u_r(z)=a(z),\;z\in G_r$.

It follows from (\ref{grad1}) that $|\<s_r,\,\zeta-z\>|$
are uniformly bounded from below by some positive
value if $z\in K,\, \zeta\in\bar G_r,\, |\xi|\ge r$ and $r\ge
2\sup\{|x|:\,x+iy\in K\}$. Hence the integral
$$
C_m\int_{G_r\cap\{\zeta:\,|\xi|>r\}}\exp\<z-\zeta,\,\zeta-z\>
{a\wedge\<s_r,\,d\zeta\>\wedge(d_\zeta\<s_r,\,d\zeta\>)^{m-1}\over
\<s_r,\,\zeta-z\>^m}.
$$
and all its derivatives in $z$ tend to zero as $r\to\infty$. For
$|\xi|<r$ we have $\rho_r(\zeta)=\rho(\zeta)$ and the function
$s_r$ coincides with $s=\overline{\<\grad_\zeta
\rho(\zeta),\,\zeta-z\>}
\grad_\zeta\rho(\zeta)+|\rho(\zeta)|(\bar\zeta-\bar z)$. Therefore
the functions $u_r$ tend to the integral \be\label{BA2} A
a(z)=C_m\int_{T_\O}\exp\<z-\zeta,\,\zeta-z\>
{a\wedge\<s,\,d\zeta\>\wedge(d_\zeta\<s,\,d\zeta\>)^{m-1}\over
\<s,\,\zeta-z\>^m}. \en as $r\to\infty$ uniformly in $z\in K$, and
the function $Aa(z)$ belongs to the class $C^1$ and $\bp
Aa(z)=a(z)$ for all $z\in T_\O$.

It is clear that the assertion a) is fulfilled. Since the function $\rho(\zeta)$
does not depend on $\xi$, the kernel in
the integral (\ref{BA2}) depends only on the variables $y,\,\eta,
\,x-\xi$, and the substitution $\zeta\mapsto\zeta-t$ for all
$t\in\R^m$ gives the assertion b). Finally, it follows from
(\ref{grad3}) that
$$
|\<s(z,\zeta),\,\zeta-z\>|\ge
\<\grad_\eta\kappa(\eta),\eta-y\>\ge C_1>0
$$
for $z\in T_{\O'},\, \zeta\in T_\O,\,
\dist(\eta,\,\partial\O)<\d$,
and
$$
|\<s(z,\zeta),\,\zeta-z\>|\ge
|\kappa(\eta)||\zeta-z|\>\ge C_2|\zeta-z|
$$
for $z\in T_{\O'},\, \zeta\in T_\O,\,\dist(\eta,\,\partial\O)\ge\d$.
Moreover, for $z,\,\zeta\in T_\O$ all the coefficients of the form
$\<s(z,\zeta)d\zeta\>$ have the
upper bound $C_3|\zeta-z|$, and the coefficients of the form
$d_\zeta\<s(z,\zeta),\,\zeta-z\>$ have the upper bound $C_4(1+|\zeta-z|)$.
Hence the integrand in (\ref{BA2}) estimates from above as
$$
C_5\exp(-|x-\xi|^2)\sup_{\zeta\in T_\O}|a(\zeta)|
{|\zeta-z|(1+|\zeta-z|)^{m-1}\over\min\{1,\,|\zeta-z|^m\}}
$$
for all $z\in T_{\O'}$, and we obtain the assertion c).

In the general case, it is sufficient to take an arbitrary convex domain
$\tilde\O$, $\O'\Subset\tilde\O\Subset\O$, such that the operator
$A$ with the desired properties for the pair ($\O',\,\tilde\O$) has been
already constructed. \bs

\medskip

Our main result of this section is as follows.
\begin{th} \label{t-bund2}
Let $\F$ be a line bundle over $\hRO$ with holomorphic transition
functions $\{g_{\a,\beta}\}$. If for some $y_0\in\O$ there exists
a section of $\F$ over $\hR\times\{y_0\}$ which does not vanish on
$\hR\times y_0$, then $\F$ has a global holomorphic section
vanishing nowhere.
\end{th}

As an immediate consequence of theorems \ref{t-bund} and
\ref{t-bund2} we get

\begin{cor}
\label{t-rest}(for $m=1$ see \cite{FRR1}). If the restriction of a divisor $d\in\AM$ to the
domain $T_{\tilde\O}\subset T_\O$ is the divisor of a function
$g\in APH(T_{\tilde\O})$ (for example, the projection of $\supp
\,d$ on $\O$  is not dense in $\O$), then $d$ is the divisor of a
function $\fHT$.
\end{cor}

\ws Since the domain $\O$ is contractible to a point $\{y_0\}$,
the bundle $\F$ is isomorphic to the restriction of $\F$ to the
base $\hR\times\{y_0\}$ (see \cite{H2}, Chapter 3, Theorem 4.7).
Therefore there exists some global section of $\F$ vanishing
nowhere, and the Chern class $c(\F)\in H^2(\hRO,\Z)$ is zero.
Hence the cocycle $\{c_{\a,\beta,\g}\}$ from (\ref{coh}) is a
coboudary for a sufficiently fine covering $\{\o_\a\}$ of $\hRO$,
i.e., there exists a 1-cochain $\{b_{\a,\beta}\}$ with the
property
$$
c_{\a,\beta,\g}=b_{\a,\beta}+b_{\beta,\g}+b_{\g,\a}.
$$
Moreover,
$$
b_{\a,\beta}=-b_{\beta,\a}.
$$
Taking into account (\ref{coh}), we see that the functions
\be \label{tran1}
F_{\a,\beta}(\tau,y)=(1/2\pi i)\log g_{\a,\beta}(\tau,y)
\en
satisfy the equalities
\be \label{tran2}
F_{\a,\beta}(\tau,y)=-F_{\beta,\a}(\tau,y),\qquad
(\tau,y)\in\o_\a\cap\o_\beta,
\en
and
\be \label{tran4}
F_{\a,\beta}(\tau,y)+F_{\beta,\gamma}(\tau,y)+F_{\gamma,\a}(\tau,y)=0
\qquad (\tau,y)\in\o_\a\cap\o_\beta\cap\o_\gamma,
\en
for a suitable choice of the branches of the logarithm.

Fix a convex domain $\O'\Subset\O$. We can assume that the covering
$\{\o_\a\}$ of the set $\hR\times\overline{\O'}$ has the form
$\o_\a=U_n\times B_k$, where $\{U_n\}_{n=1}^N$ is an open covering
of $\hR$, and $\{B_k\}_{k=1}^K$ is an open covering of
$\overline{\O'}$. We can also assume that the
oscillations of all the functions $F_{\a,\beta}(\tau,y)$ on each
$\o_\a$ to be less than $1/4$.

Let $\{\chi_n(\tau)\}_{n=1}^N$ be a partition of unity corresponding to
the covering $\{U_n\}_{n=1}^N$, $\{\xi_k(y)\}_{k=0}^K$ be a partition of
unity corresponding to the covering $\{B_k\}_{k=0}^K$. Further, suppose
$\nu(s)\ge 0$ is a $C^\infty$-smooth function with the support in the ball
$|s|<\e$ and such that $\int_{\R^m}\nu(s)ds=1$. Since the covering
$\{U_n\}$ is finite, we see that for $\e$ sufficiently small the functions
$$
\tilde\chi_n(\tau)=\int_{\R^m}\chi_n(\tau\cdot js)\nu(s)ds
$$
form a partion of unity subjected to the covering $\{U_n\}$, too.
Put for $(\tau,y)\in\o_\a$
\be \label{unit}
F_\a(\tau,y)=\sum_{l,k}\tilde\chi_l(\tau)\xi_k(y)F_{(l,k),\a}(\tau,y).
\en
It follows from (\ref{tran2}) and (\ref{tran4}) that for
$(\tau,y)\in\o_\a\cap\o_\beta$,
\be \label{tran3}
F_\a(\tau,y)-F_\beta(\tau,y)=F_{\a,\beta}(\tau,y).
\en
In addition, all the functions $F_\a(\tau\cdot jx,y)$ are
$C^\infty$-smooth on their domains of definition. Put
$$
a(\tau,x,y)=\bar\partial F_\a(\tau\cdot jx,y).
$$
for $(\tau\cdot jx,y)\in\o_\a$ and fixed $\tau$.
Since the function $F_{\a,\beta}(\tau,y)$ is holomorphic, relation
(\ref{tran3}) implies that all the coefficients of
 $a(\tau,x,y)$ are well-defined and uniformly bounded on
$\hR\times\R^m\times\O'$. Besides, $a(\tau\cdot
jx,0,y)=a(\tau,x,y)$ and $\bp_za=0$. Therefore it follows from
Propositions \ref{t-d} that for each convex domain
$\O''\Subset\O'$ there exists a function $b(\tau,x,y)$, uniformly
continuous on $\hR\times\R^m\times\O''$ and such that $\bp_zb=a$
and $b(\tau,x,y)=S_xb(\tau,0,y)=b(\tau\cdot jx,0,y)$. In view of
(\ref{tran3}), the functions
$$
\p_\a(\tau,y)=\exp 2\pi i[F_\a(\tau,y)-b(\tau,0,y)],
\ (\tau,y)\in\o_\a,
$$
satisfy the equations
$$
\p_\a(\tau,y)=g_{\a,\beta}(\tau,y) \p_\beta(\tau,y),\quad
(\tau,y)\in\o_a\cap\o_\beta,
$$
hence they form a holomorphic, vanishing nowhere section of $\F$
over $\hR\times\O''$.

Let $\O_n, \ n=1,2,\dots$, be convex domains such that
$\O_n\Subset\O_{n+1}$, $\cup_n\O_n=\O$, and $\{\o_\a\}$ be a
sufficiently fine covering of $\hRO$. We have just proved that for
each $n$ there exists a holomorphic, vanishing nowhere section
$\{\p_\a^n(\tau,y)\}$ over $\hR\times\O_n$. It is obvious that the
functions
$$
\bar h_n(\tau,y)={\p_\a^{n+1}(\tau,y)\over \p_\a^n(\tau,y)},\quad
(\tau,y)\in\o_\a,
$$
are well-defined, holomorphic, and vanishing nowhere on the set
$\hR\times\O_n$. It follows from Proposition \ref{t-ex1} that
$h_n(z)=\bar h_n(jx,y)\in APH(T_{\O_n})$ and $h_n(z)\neq 0$ on $T_{\O_n}$.
Therefore (for $m=1$ see \cite{C}, for $m>1$ see \cite{F1})
\be
\label{mot}
h_n(z)=\exp[i\<z,c_n\>+\kappa_n(z)]
\en
with $c_n\in\R^m,\ \kappa_n\in APH(T_{\O_n})$. Take an exponential sum
$P_n(z)$ of the type (\ref{sum2}) such that
$|\kappa_n(z)-P_n(z)|<2^{-n}$
on $T_{\O_{n-1}}$. It is easy to prove that the function
$$
f_n(z)=\prod_{k=1}^{n-1}\exp[-i\<z,c_k\>-P_k(z)]
\prod_{k=n}^\infty h_k(z)\exp[-i\<z,c_k\>-P_k(z)]
$$
belongs to the class $APH(T_{\O_{n-1}})$. Besides,
$f_n(z)/f_{n+1}(z)=h_n(z)$.
Therefore the function
$\Psi_\a(\tau,y)=\p_\a^n(\tau,y) \bar f_n(\tau,y)$
is well-defined on $\o_\a$ whenever $\o_\a\subset\O_{n-1}$.
Hence the section $\{\Psi_\a(\tau,y)\}$  is well-defined,
holomorphic, and vanishing nowhere on the set $\hRO$. \bs

\medskip

{\bf Remark 1.} Instead of the existence of a section
over $\hR\times y_0$ vanishing nowhere we may assume that
(\ref{tran2}) and (\ref{tran4}) is true for $y=y_0$.

\medskip

{\bf Remark 2.} The condition $h_n\in APH(T_{\O_n},\G)$ implies
$k_n\in APH(T_{\O_n},\G)$ and $c_n\in\Gamma$ (see \cite{F1}). Therefore, since
the class $APH(T_\O,\G)$ is closed with respect to the uniform convergence,
the previous theorem is true for line bundles over
$\hat\G\times\O$.

\bigskip

\section{Chern class of almost periodic divisors}
\label{sec4}

\bigskip

As above, take a covering $\{U_\a\times B_k\}$ of the set $\hRO$
such that $\{U_\a\}_{\a\in A}$ is an open covering of $\hR$ and
$\{B_s\}_{0\le s<\infty}$ is an open covering of $\O$. Suppose the
functions $r_{(\a,s)}(\tau,y)$ define a divisor $d\in\AM$ on the
sets $\{U_\a\times B_s\}$; then the functions
$g_{(\a,s),(\beta,n)}(\tau,y)=r_{(\a,s)}(\tau,y)/r_{(\beta,n)}(\tau,y)$
are transition functions for the line bundle $\F_d$ over $\hRO$.
Fix $y_0\in B_0$; the function $g_{(\a,0),(\beta,0)}(\tau,y_0)$
defines the line bundle $\F_{d,y_0}$ over $\hR$. If we take
another functions $\tilde r_{(\a,s)}(\tau,y)$ defining the divisor
$d$ over $U_\a\times B_s$, then $\tilde r_{(\a,s(}(\tau,y)=
h_{(\a,s)}(\tau,y) r_{(\a,s)}(\tau,y)$ with some functions
$h_\a(\tau,y)$ vanishing nowhere. Therefore a unique Chern class
is assigned to the divisor $d$.

We will say that the Chern class $c(\F_{d,y_0})$ is the {\sl Chern
class $c(d)$ of the divisor $d$}.

Note that the sum of divisors $d_1,\,d_2$ is defined by the
product of the functions defining the divisors $d_1,\,d_2$; hence,
$c(d_1+d_2)=c(d_1)+c(d_2)$.

Also, if $d$ is the divisor of a function $\fHT$, then by Theorem
\ref{t-bund1} the bundle $\F_{d,y_0}$ has a section
$\p_\a(\tau)$ vanishing nowhere. Take in (\ref{coh})  \ $\log
g_{\a,\beta}(\tau)=\log\p_\a(\tau)-\log\p_\beta(\tau)$
for all $\a,\beta$, then we obtain $c(d)=0$. On the other hand,
if $c(d)=0$, then for a sufficiently fine covering of $\hR$ and
for a suitable branches of logarithm, the left-hand side in
(\ref{coh}) is zero for all $\a,\beta,\g$. Now, using Theorem
\ref{t-bund} and Remark 1 to Theorem \ref{t-bund2}, we come to the
following result.

\begin{th}
\label{t-Cher}
$d\in\AM$ is the divisor of a function $\fHT$ iff $c(d)=0$.
\end{th}

Suppose that a domain $\O\subset\R^m$ is stable under the map
$L:\,y\mapsto 2y_0-y$. If $\tilde d$ is the image of a divisor
$d\in\AM$ under this map, then the holomorphic functions $\tilde
r_\a(\tau,y)=\overline{r_\a(\tau,2y_0-y)}$ define this divisor on the
set $U_\a\times B_0$ (we may assume that $L(B_0)=B_0$). Then we have
$\log\tilde g_{(\a,0)(\beta,0)}(\tau,y_0)=
-\log g_{(\a,0)(\beta,0)}(\tau,y_0)$.
This implies
$c(\tilde d)=-c(d)$. Using Theorem \ref{t-Cher} and Corollary
\ref{t-rest}, we get

\begin{th}
\label{t-sym}
If the restriction of a divisor $d\in\AM$ to the ball $B\subset\O$ with
center in $y_0$ is stable under the map $y\mapsto 2y_0-y$, then $d$ is the
divisor of a function $\fHT$.
\end{th}

Suppose $\sp d\subset\G$, where $\G$ is a subgroup of $\R^m$. As
above, we can introduce the Chern class $c_\G(d)\in H^2(\hG,\Z)$
and prove the following statement.

\begin{th}
\label{t-CherG}
$d\in APM_{1,1}(T_\O,\G)$ is the divisor of a function $f\in APH(T_\O,\G)$
iff $c_\G(d)=0$.
\end{th}

Consider a relation between $c_\G(d)$ and $c(d)$. Let $\iota$ be
the identity mapping of $\R^m$ with the topology $\beta$ to $\R^m$
with the topology $\beta(\G)$. Extend $\iota$ to a continuous
mapping $\bar\iota$ of $\hR$ to $\hG$. It is easy to prove that
$\bar\iota$ is an open surjection (actually, $\bar\iota$ is an
epimorphism of the group $\hR$ to $\hG$). Let $\bar d$ be the
continuous mapping from Proposition \ref{t-ex2}, $\bar d_\G$ be
the corresponding mapping of $\hG$ to $\AM$ (see Remark after
Proposition \ref{t-ex2}). It can be easily checked that $\bar
d=\bar d_\G\circ\bar\iota$. Hence \be \label{mono}
c(d)=\iota^*(c(d)), \en where the homomorphism
$\iota^*:\,H^2(\hG,\Z)\mapsto H^2(\hR,\Z)$ is the lifting of
$\bar\iota$. Since $\bar\iota$ is an open surjection, $\iota^*$ is
a monomorphism. Using Theorem \ref{t-Cher}, Theorem \ref{t-CherG}
and (\ref{mono}), we derive the following result.

\begin{th}
\label{t-1}
If $d$ is the divisor of a function $f\in APH(T_\O)$ and $\sp d\subset\G$,
then $d$ is the divisor of a function $f_1\in APH(T_\O,\G)$
\end{th}

It follows from the next theorem that the case
$\G=\Lin_\Z\{\l_1,\dots,\l_N\}$ with the vectors
$\l_1,\dots,\l_N\in\R^m$ linearly independent over $\Z$ is most
important.

\begin{th}
\label{t-Cher1} For each $d\in\AM$ there exists a divisor
$d'\in\AM$ with the spectrum in the group
$\Lin_\Z\{\l_1,\dots,\l_n\}$ such that the both divisors have the
same Chern class.
\end{th}

\ws Let $\{U_\a\}$ be any sufficiently fine covering of $\hR$. We
may assume that
$$
U_\a=\{\tau\in\hR:\,
(\tau_{\l_1},\dots,\tau_{\l_n})\in\tilde U_\a\},
$$
where  $\l_1,\dots,\l_n\in\R^m$ are linearly independent over $\Z$,
the number $n$ and the vectors $\l_1,\dots,\l_n$ are the same for
all $\a$, $\tau_{\l_1},\dots,\tau_{\l_n}$ are the coordinates of
$\tau\in\hR\subset\prod_{\l\in\R^m}\T_\l$,
$\tilde U_\a$ is an open set in the torus $\T^n$. Clearly, the group
$H^2(\{U_\a\},\Z)$ is isomorphic to the group $H^2(\{\tilde
U_\a\},\Z)$. Passing to the inductive limit with respect to the refinement of
coverings, we obtain a monomorphism
$$
\iota^*:\,H^2(\hG,\Z)\mapsto
H^2(\hR,\Z)
$$
with $\G=\Lin_\Z\{\l_1,\dots,\l_n\}$, and the equality
$$
H^2(\hR,\Z)=\cup_\G\iota^*H^2(\hG,\Z),
$$
where $\G$ runs over all
the subgroups of the type $\Lin_\Z\{\l_1,\dots,\l_n\}$. It will
be proved in the next section that for any element $c\in
H^2(\hG,\Z)$ with $\G=\Lin_\Z\{\l_1,\dots,\l_n\}$ there exists a
divisor $d'\in\AM$ with the Chern class $c$. \bs

\medskip

Let $d\in APM_{1,1}(T_\O,\G)$ with
$\G=\Lin_\Z\{\l_1,\dots,\l_N\}$, the vectors
$\l_1,\dots,\l_N\in\R^m$ being linearly independent over $\Z$.
Here $\hG=\T^N=\{\zeta\in\C^N:\,\zeta^l=e^{2\pi iu^l},\
l=1,\dots,N\}$. By definition, put
\be \label{map}
P:\,u\in\R^N\mapsto (e^{2\pi iu^1},\dots,e^{2\pi iu^N}).
\en
Let $\{U_\a\}$ be any sufficiently fine open covering of $\T^N$
with connected $U_\a$, $\F$ be a
line bundle over $\T^N$ with transition functions
$\{g_{\a,\beta}(\zeta)\}$. Take an arbitrary connected component
$\o_{\a,0}$ of the set $P^{-1}(U_\a)$ and put
$$
\o_{\a,k}=\o_{\a,0}+k,\quad k\in\Z^N.
$$
It can be assumed that
$\o_{\a,k}\cap\o_{\a,k'}=\emptyset$ for all $\a$ and $k\neq k'$.
The functions $G_{(\a,k),(\beta,n)}(u)=g_{\a,\beta}(P(u)), \
u\in\o_{\a,k}\cap\o_{\beta,n}$, are transition functions for some
line bundle $\bar\F$ over $\R^N$. Any bundle over $\R^N$ has a
global section vanishing nowhere, and let $\{\Phi_{\a,k}(u)\}$ be
just the same section for $\bar\F$. Then we have
\be \label{re1}
\Phi_{\a,k}(u)=g_{\a,\beta}(P(u))\Phi_{\beta,n}(u),\quad
u\in\o_{\a,k}\cap\o_{\beta,n}.
\en
Further, by $e_1,\dots,e_N$ denote the basis vectors in $\R^N$. Since
$g_{\a,\beta}(P(u+e_l))=g_{\a,\beta}(P(u))$, we see that the functions
\be \label{re2}
\Psi_l(u)=\Phi_{\a,k+e_l}(u+e_l)/\Phi_{\a,k}(u),\quad l=1,\dots,N,
\ u\in\o_{\a,k},
\en
are well-defined on $\R^N$. Then the numbers
\be \label{entr}
m_{p,q}={1\over 2\pi i}[\log\Psi_p(u+e_q)-\log\Psi_p(u)-
\log\Psi_q(u+e_p)+\log\Psi_q(u)],\quad 1\le p,q\le N,
\en
form a skew-symmetric matrix $M$ with integral entries. This matrix
does not depend on the branches of $\log\Psi_l(u),\ 1\le l\le N$.
Moreover, since any other section $\{\tilde\Phi_{\a,k}(u)\}$
vanishing nowhere on $\R^N$ has the form
$\tilde\Phi_{\a,k}(u)=\Xi(u)\Phi_{\a,k}(u)$ with $\Xi(u)\neq 0$
on $\R^N$, we see that $M$ is defined only by the bundle $\F$. It
is easy to prove that if matrices $M$ and $\tilde M$ correspond
to the bundles $\F$ and $\tilde\F$ over $\T^N$ with transition
functions $\{g_{\a,\beta}\}$ and $\{\tilde g_{\a,\beta}\}$, then
the matrix $M-\tilde M$ corresponds to the bundle with the
transition functions $\{g_{\a,\beta}/\tilde g_{\a,\beta}\}$. If
the bundle $\F$ has a global section $\{\p_\a(\tau)\}$ vanishing
nowhere, then the functions $\Phi_{\a,k}(u)=\p_\a(P(u)),
u\in\o_{\a,k}$ form a global section of $\tilde\F$ over $\R^N$.
Since this section satisfies the conditions
\be \label{cond}
\Phi_{\a,k+e_l}(u+e_l)=\Phi_{\a,k}(u),\quad u\in\o_{\a,k}
\en
for all $\a,\ k\in\Z^N,\,l=1,\dots,N$,  we see that $\Psi_l(u)\equiv
1,\ l=1,\dots,N$, and $M=0$. In particular, if the bundles $\F$
and $\tilde\F$ over $\T^N$ have the same Chern class, then the
same matrix $M$ corresponds to these bundles. Thus the mapping
$c\mapsto M$ is a well-defined homomorphism of $H^2(\T^N,\Z)$ to
the additive group of all skew-symmetric matrices with integral
entries.

Let us check that this mapping is injective.

Suppose that the matrix $M=0$ corresponds to the bundle $\tilde\F$
over $\R^N$ and
$\{\Phi_{\a,k}(u)\}$ is a global section of $\tilde\F$ vanishing nowhere.
Put
\begin{eqnarray*}
H_1(u)=\left\{
\begin{array}{ll}
-\sum_{n=1}^{[u^1]}\log\Psi_1(u-ne_1)-(u^1-[u^1])\log\Psi_1(0,'u)
\quad &{\rm for}\; u^1\ge 1,\\
-u^1\log\Psi_1(0,'u)\quad &{\rm for}\; 0\le u^1<1,\\
\sum_{n=o}^{-[u^1]-1}\log\Psi_1(u-ne_1)-(u^1-[u^1])\log\Psi_1(0,'u)
\quad &{\rm for}\;u^1<0,
\end{array}
\right.
\end{eqnarray*}
and for all $\a,\,k$
$$
\tilde\Phi_{\a,k}(u)=\Phi_{\a,k}(u)\exp H_1(u),
$$
where  $u=(u^1,'u)$, $u^1\in\R,\ 'u\in\R^{N-1}$, and $[u^1]$ is
the integral part of $u^1$. It is easy to see that the function
$H_1(u)$ is continuous and satisfies the equality
$H_1(u+e_1)-H_1(u)=-\log\Psi_1(u)$. It follows from (\ref{re2})
that the section $\{\tilde\Phi_{\a,k}(u)\}$ satisfies (\ref{cond})
for $l=1$.

Further, suppose that $\{\Phi_{\a,k}\}$ satisfies
(\ref{cond}) for some $l'\neq 1$. By (\ref{re2}), we can take
$\log\Psi_{l'}(u)\equiv 0$. Using the equality $M_{1,l'}=0$, we have
$\log\Psi_1(u+e_{l'})-\log\Psi_1(u)=\log\Psi_{l'}(u+e_1)-\log\Psi_{l'}(u)=0$
and $H_1(u+e_{l'})-H_1(u)=0$. Hence the
section $\Phi_{\a,k}\exp H_1(u)$ satisfies (\ref{cond}) with
$l=1,\,l=l'$. Therefore arguing as above, we can "improve"
the section in each coordinate sequentially and obtain a section
$\{\tilde\Phi_{\a,k}\}$
satisfying (\ref{cond}) for all $l$. Now it follows that the
functions $\p_\a(\zeta)=\hat\Phi_{\a,k}(P^{-1}(\zeta))$ are
well-defined on $U_\a\subset\T^N$ for all $\a$ and form a global
section vanishing nowhere.

Below we will prove that each skew-symmetric matrix with integer
entries corresponds to some Chern class of a purely periodic
divisor, therefore the constructed mapping $c\mapsto M$ is an
isomorphism. Thus we will identify a Chern class $c\in
H^2(\T^N,\Z)$ with the matrix $M$.

\bigskip

\section{Periodic divisors and Completing Theorem}
\label{sec5}

Let $D$ be a divisor in the domain $T_G=\R^N+iG$, where
$G\subset\R^N$ is a convex domain. The divisor is called {\sl
$N$-periodic} if there exist $N$ vectors $u_l\in\R^N$, linearly
independent over $\R$ and such that $D(w+u_l)=D(w), \
l=1,\dots,N$. We can assume without loss of generality that these
vectors are $e_1,\dots,e_N$. Then $\sp D\subset 2\pi\Z^N$. Remark
that we do not exclude the case of divisors depending on $k<N$ coordinates.

Let $F(w)\in H(T_\O)$ be an arbitrary function with the divisor
$D$. The mapping $\bar D(\zeta)=\dd\log|F(w+{\log\zeta\over 2\pi
i})|$, where $\log\zeta=(\log\zeta^1, \dots,\log\zeta^N)$, is
well-defined; by Proposition \ref{t-con}, it continuously maps
from $\T^N$ to $M(T_G)$.

Let $\L$
%$\L=\left(\begin{array}{c}\l_1\\ \dots\\\l_N\end{array}\right)$
be the matrix with rows
$\l_1,\dots,\l_N\in\R^m$ linearly independent over $\Z$, $\O$ be
the set $\{y\in\R^m:\,\L y\in G\}$. If
$\G=\Lin_\Z\{\l_1,\dots,\l_N\}$, then the map $j:\,\R^m\mapsto\hG
(=\T^N)$ is $P\circ\L$, where $P$ is defined in (\ref{map}). The
mapping $\bar d(\zeta)=\dd\log|F(\L z+{\log\zeta\over 2\pi i})|$
continuously takes $\T^N$ to $M_{1,1}(T_\O)$ too, and $$\bar
d(\zeta\cdot jt) =\bar d(\zeta\cdot P(\L t))=\dd\log|F(\L z+\L t+
{\log\zeta\over 2\pi i})|=S_t\bar d(\zeta)$$ for all $t\in\R^m$.
Using Proposition \ref{t-ex2}, we see that the divisor $d=\bar
d(1)$ of the function $F(\L z)$ is almost periodic with spectrum
in $\G$.

Take a sufficiently fine covering $\{U_\a\}_{\a\in A}$ of $\T^N$.
As in the previous section, define open sets
$\o_{\a,k}\subset\R^N,\ k\in\Z^N$. The functions
$r_\a(\zeta,y)=F(u+i\L y)$ for $u\in\o_{\a,k}, P(u)=\zeta,\
k=k(\a)$, define the divisor $d$ on $U_\a\times\O$ for any mapping
$k(\a)$ from $A$ to $Z^N$. Then the functions
$g_{\a,\beta}(\zeta,y)=r_\a(\zeta,y)/r_\beta(\zeta,y)$ are
transition functions for the line bundle $\F$ over $\T^N\times\O$.

Fix $y_0\in\O$. Then the functions
$\Phi_{\a,k}(u)=r_\a(P(u),y_0)/F(u+i\L y_0)$ form a global
holomorphic section of some line bundle over $\R^N$, vanishing
nowhere. Now the equalities (\ref{re2}) and (\ref{entr}) define
the functions $\Psi_l(u),\ l=1,\dots,N$ and the matrix $M$. Since
$P(u+e_l)=P(u)$ for all $l=1,\dots,N$, we have
$$
\Psi_l(u)=F(u+i\L y)/F(u+e_l+i\L y),\quad l=1,\dots,N.
$$
The entries of the matrix $M$ are integer and depend continuously
on $\Psi$,  therefore we can change $\Psi_l(u)$ to the functions
\be \label{inc} \tilde\Psi_l(w)=F(w)/F(w+e_l) \en for all fixed
$w\in G$.

Since all projections of the set $\{w:\,F(w)=0\}$ to the
hyperplanes $\Im w^l=0,\,l=1,\dots,N$, have zero Lebesque measure,
we see that the rectangle $\Pi_{p,q}(w)$ with vertices
$w,\,w+e_p,\,w+e_p+e_q,\,w+e_q$ does not intersect $\supp D$ for
a.a. $w\in T_G$. Therefore for a.a. $w'$, \be
\begin{array}{l@{=}c@{=}r}
-\log\tilde\Psi_p(w')&\log F(w'+e_p)-\log F(w')&
\Delta_{[w',w'+e_p]}\log F(w),\\
-\log\tilde\Psi_p(w'+e_q)&\log F(w'+e_p+e_q)-\log F(w'+e_q)&
\Delta_{[w'+e_q,w'+e_p+e_q]}\log F(w),
\end{array}
\en
where $\Delta_{[a,b]}f$ means the increment of the function $f$ on the
segment $[a\,b]$. Similarly,
\be
\begin{array}{l@{=}r}
-\log\tilde\Psi_q(w')&\Delta_{[w',w'+e_q]}\log F(w),\\
-\log\tilde\Psi_q(w'+e_p)&\Delta_{[w'+e_p,w'+e_p+e_q]}\log F(w),
\end{array}
\en
Hence
$$
m_{p,q}=(1/2\pi)\Delta_{\Pi_{p,q}(w)}\Arg F(w).
$$

{\bf Remark.} In \cite{R4} L.Ronkin proves that a necessary and
sufficient condition for a periodic divisor to be the divisor of a
holomorphic periodic function is $M=0$; in \cite{R6} he shows that
a section of a periodic divisor $D$ (we denote this section by
$d$) is the divisor of a holomorphic almost periodic function iff
$D$ is the divisor of a holomorphic periodic function. Here we
have shown that the Chern classes of $d$ and $D$ coincide, hence
Ronkin's result follows from Theorem \ref{t-Cher}.

\medskip

{\bf Example.} Let $\phi(\zeta)$ be an entire function in the
plane $\C$ with simple zeroes at all points with integer
coordinates. Put for $w=(w^1,\dots,w^N)\in\C^N$ and fixed
$p,q\in\N,\,1\le p<q\le N$, \be \label{exam}
F_{p,q}(w)=\phi(w^p+iw^q). \en By $D_{p,q}$ denote the divisor of
$F$. Since $\tilde\Psi_l(w)\equiv 1$ for $l\neq s,t$, we see that
all the entries of the corresponding matrix $M_{p,q}$ (except for
$m_{p,q}$ and $m_{q,p}=-m_{p,q}$) vanish. Take
$w_0=\{-1/2,\dots,-1/2\}$. Since $F_{p,q}(w)\neq 0$ for
$w\in\Pi_{p,q}(w_0)$, we have \be \label{inc1}
m_{p,q}=(1/2\pi)\Delta_{\Pi_{p,q}(w_0)}\Arg F_{p,q}(w)=
(1/2\pi)\Delta_L\Arg\phi(w)=1, \en here $L$ is the rectangle in
$\C$ with vertices $-1/2-i/2,\,1/2-i/2,\,1/2+i/2,\,-1/2+i/2$.

The divisor $d_{\l_p,\l_q}$ of the function
$$
F_{p,q}(\L z)=\phi(\<z,\l_p\>+i\<z,\l_q\>),\quad z\in\C^m,
$$
has the same Chern class (as above, $\L$ is the $N\times m$-matrix
with rows $\l_1,\dots,\l_N\in\R^m$ linearly independent over
$\Z$). It follows easily that if $\l_p,\,\l_q$ are linearly
independent over $\R$, then $d_{\l_p,\l_q}$ is $N$-periodic.
Otherwise, the support of $d_{\l_p,\l_q}$ is a union of parallel
complex hyperplanes.

\medskip

Now we describe Chern classes of divisors as finite sums of
special type.

For vectors $\l,\,\mu\in\R^m$, denote by $\l\wz\mu$ the Chern
class of the divisor $d_{\l,\mu}$ of the function
$\phi(\<z,\l\>+i\<z,\mu\>)$. First let $\l,\,\mu$ be linearly
independent over $\Z$. The permutation of $\l$ and $\mu$
corresponds to the rearrangement of $p$ and $q$ in (\ref{exam}),
therefore we get $-1$ in (\ref{inc1}). Thus \be \label{rul1}
\mu\wz\l=-\l\wz\mu. \en Taking $\l\mapsto -\l$ corresponds to
$w^p\mapsto -w^p$ in (\ref{exam}), therefore we get $-1$ in
(\ref{inc1}) again. Hence \be \label{rul2} (-\l)\wz\mu=-\l\wz\mu.
\en Changing $\l$ to $k\l,\,k\in N$, corresponds to changing $w^p$
to $kw^p$, so we obtain $m_{p,q}=k$ in (\ref{inc1}). Since other
coefficients of $M$ are zero, we get \be\label{rul3}
(k\l)\wz\mu=k(\l\wz\mu); \en if we change $\l$ to $(n/k)\l$ here,
then we get (\ref{rul3}) for all rational $n/k$.

Further, let a matrix $\tilde M$ represent the Chern class of the
divisor of the function $\tilde F(w)=\phi(w^p+w^r+iw^q)$. Arguing
as above, we get $\tilde m_{p,q}=\tilde m_{r,q}=1,\ \tilde
m_{p,r}=0$. Therefore the Chern class of this divisor equals the
sum of the Chern classes of the divisors of the functions
$\phi(w^p+iw^q)$ and $\phi(w^r+iw^q)$. Let
$\l',\,\l'',\,\mu\in\R^m$ be linearly independent over $\Z$. If we
take a matrix $\L$ with the rows
$\l_p=\l',\,\l_r=\l'',\,\l_q=\mu$, we see that the Chern class of
the divisor of the function $\phi(\<z,\l'+\l''\>+i\<z,\mu\>)$ equals
the sum of the Chern classes of divisors $d_{\l',\mu}$ and
$d_{\l'',\mu}$. We have
\be \label{rul4}
(\l'+\l'')\wz\mu=\l'\wz\mu+\l''\wz\mu.
\en
Consider the divisor
of the function $\phi(kw^p+iw^p)$. Since this function depends
only on one variable, we see that the zero matrix $M$ corresponds
to the divisor $d_{k\l,\l}$. Therefore we have
\be \label{rul5}
\l\wz\mu=0 ,
\en
if the vectors $\l,\,\mu$ are linearly dependent
over $\Z$. Now, arguing as in (\ref{inc1}), we obtain that the
divisor of the function $\phi(w^p+w^q+iw^q)$ has the same matrix
$M$ as the divisor of the function $\phi(w^p+iw^q)$, hence
$$
(\l+\mu)\wz\mu=\l\wz\mu.
$$
Therefore if $\l''=r\l'+s\mu$ for rational $r$ and $s$, then
\begin{eqnarray*}
(\l''+\l')\wz\mu=s((r+1)/s\l'+\mu)\wz\mu=r\l'\wz\mu+\l'\wz\mu\\
=s(r/s\l'+\mu)\wz\mu+\l'\wz\mu=\l''\wz\mu+\l'\wz\mu.
\end{eqnarray*}
Thus (\ref{rul4}) is true for all $\l',\,\l'',\,\mu$.

Since $M=\sum_{p<q}m_{p,q}M_{p,q}$, we get that the Chern class
of each divisor with spectrum in
$\Lin_\Z\{\l_1,\dots,\l_N\}$ has the form
\be \label{Cher}
\sum_{1\le p<q\le N}m_{p,q}\l_p\wz\l_q.
\en

Actually, this statement realizes a well-known isomorphism between
$H^2(\hR,\Z)$ and $\R^m\wz\R^m$ (see \cite{H}).

\medskip

Now we can prove that every almost periodic divisor is
complemented to the divisor of some holomorphic almost periodic
function.

\begin{th}
\label{t-comp} The Chern class of each divisor $d\in\AM$ can be
represented as a finite sum $\sum\l_s\wz\mu_s$; the divisor
$d+\sum d_{\mu_s,\l_s}$ is the divisor of some $\fHT$. If $m>1$,
then we can choose $\mu_s,\,\l_s$ in such a way that all the
divisors $d_{\mu_s,\l_s}$ are periodic.
\end{th}

\ws The first statement follows from (\ref{Cher}), (\ref{rul3}),
and Theorem \ref{t-Cher1}. Using
(\ref{rul1}) and Theorem \ref{t-Cher}, we get the second
statement of our theorem. Finally, note that for $m>1$ one can
take a basis in $\G'\supset\G=\Lin_\Z\{\l_1,\dots,\l_N\}$ such
that any two elements of the basis are linearly independent over $\R$.
Namely, choose
$$
\l_0\not\in\cup_{1\le j,k\le N}\{t\l_j+(1-t)\l_k:\,t\in\R\}\cup
\cup_{1\le j\le N}\{t\l_j:\,t\in\R\}
$$
and take vectors $\l_0,\,\l_1-\l_0,\dots,\l_N-\l_0$
as a basis of the group $\G'=\Lin_\Z\{\l_0,\l_1,\dots,\l_N\}$.\bs

\bigskip

\section{Almost periodic mappings into projective space}
\label{sec6}

\bigskip

Let $F(z)$ be a holomorphic mappings of $T_\O$ into the
projective space $\CP$. Using the homogeneous coordinates, it can
be written in the form $[f^0(z):f^1(z):\dots:f^k(z)]$ with
holomorphic functions $f^l$ on $T_\O$ without common zeroes.
These functions is well-defined up to common holomorphic factor
vanishing nowhere on $T_\O$. We will say that divisors $d_l$ of
the functions $f^l$, $l=0,\dots,k$, are coordinate divisors for $F$
(we suppose that the mapping $F$ is not degenerate, i.e., each
coordinate $f^l$ is not identically zero). We assume that $\CP$ is
equipped with the Fubini-Study metric.

In the case $k=1$ we can interpret $F(z)$ as a meromorphic
function on $T_\O$ with disjoint sets of zeroes and poles  and
the spherical metric on the set of values of $F$.

We say that a mapping $F:\ T_\O\mapsto\CP$ is almost periodic if
the collection of the shifts $\{S_tF\}_{t\in\R^m}$ is a relatively compact set
with respect to the topology of uniform convergence on every tube domain
$T_{\O'}$, $\O'\Subset\O$.

For example, $F(z)=[1:f^1(z):f^2(z):\dots:f^k(z)]$ with $f^l\in
APH(T_\O)$ is an almost periodic mapping, because the Fubini-Study
metric is equivalent to the Euclidean metric for the local chart
$W^0=1$.

Like almost periodic functions, almost periodic mappings extend to
continuous mappings from $\hRO$ to $\CP$. We will also use the
following statement.

\begin{prop}
\label{t-map1} To each holomorphic almost periodic mapping
$F:\,T_\O\mapsto\CP$ one can assign a line bundle $\F$ over $\hRO$
and $k+1$ holomorphic sections of this bundle without common zeroes,
which define coordinate divisors for $F$.
Conversely, any $k+1$ holomorphic sections of a line bundle
without common zeroes correspond to coordinate divisors of some
holomorphic almost periodic mapping from $T_\O$ to $\CP$.
\end{prop}

\ws Let $F:\,T_\O\mapsto\CP$ be a holomorphic almost periodic mapping.
Consider the function
$$
\kappa(t,y)=\sup_{x\in\R^m}\rho(F(x+iy+t),F(x+iy)),\quad t\in\R^m,\ y\in\O,
$$
where $\rho$ is the Fubini-Study distance on $\CP$. It follows
from the inequality
$$
|\kappa(t,y)-\kappa(t',y)|\le\sup_{x\in\R^m}\rho(F(x+iy+t),F(x+iy+t')),
$$
that $\kappa(x,y)\in AP(T_\O)$, hence $\kappa$ extends to a continuous
function on $\hRO$. This means that both $\kappa$ and $F$ are
uniformly continuous in the topology $\beta$. Therefore $F(z)$
extends to a continuous mapping $\bar
F(\tau,y):\,\hRO\mapsto\CP$. Fix a point $(\tau_0,y_0)\in\hRO$.
Then the values of the mapping $\bar F$ for $\tau, y$ in some
neighborhood $\o$ of this point belong to one of the local
charts, for example to $W^0=1$. Hence $\bar
F(\tau,y)=[1:\p^1(\tau,y):\dots:\p^k(\tau,y)]$ for $(\tau,y)\in\o$
with continuous functions $\p^1,\dots,\p^k$. Moreover, since
$\bar F(jx,y)=F(z)=[1:f^1(z):\dots:f^k(z)]$ for $(jx,y)\in\o$, the
functions $\p^l(\tau,y),\ l=1,\dots,k$, are holomorphic on $\o$.

Let $\{\o_\a\}$ be a sufficiently fine open covering of $\hRO$. It can be assumed
that for each $\o_\a$ at least one of the homogeneous coordinates of $\bar
F=[\p_\a^0:\p_\a^1:\dots:\p_\a^k]$ equals identically $1$ on $\o_\a$. Taking
into account Proposition \ref{t-quot}, we obtain that the
functions $g_{\a,\beta}=\p_\a^0/\p_\beta^0=\dots=\p_\a^k/\p_\beta^k$ have
no zeroes on $\o_\a\cap\o_\beta$. It is clear that these
functions define a line bundle $\F$ over $\hRO$, and functions
$\{\p_\a^l(\tau,y)\}$ for all fixed $l=0,\dots,k$, define global
holomorphic sections of $\F$. Proposition \ref{t-bund1} shows that
the coordinate divisors of $F(z)$ are defined by just these
sections.

Conversely, any holomorphic global sections
$\p^l=\{\p_\a^l(\tau,y)\},\ l=0,\dots,k$, of the line bundle $\F$
over $\hRO$ without common zeroes are homogeneous coordinates for
the continuous holomorphic mapping
$\bar F=[\p_\a^0,\dots,\p_\a^k]: \o_\a\mapsto\CP$
for each $\a$. Since
$$
\p_\a^0/\p_\beta^0=\dots=\p_\a^k/\p_\beta^k,
$$
the mapping $\bar F$ is well-defined on $\hRO$; the continuity
of $\bar F(\tau,y)$ implies almost periodicity of $F(z)=\bar F(jx,y)$. \bs

\medskip

The following theorem is the main result of this section.

\begin{th}
\label{t-map2}

In order that divisors $d_0,\dots,d_k$ be the coordinate divisors
of some holomorphic almost periodic mapping $F:\,T_\O\mapsto\CP$, it
is sufficient and necessary that the following conditions are
fulfilled:

a) $d_0,\dots,d_k$ are almost periodic divisors,

b) $d_0,\dots,d_k$ have the same Chern class,

c) for each $\O'\Subset\O$ there exists $\d>0$ such that every ball
$B(z_0,\d),\ z_0\in\O'$, intersects at most $k$ supports of the
divisors $d_0,\dots,d_k$.
\end{th}

\ws Let $F:\,T_\O\mapsto\CP$ be a holomorphic almost periodic mapping.
It follows from Proposition \ref{t-map1} that there exist holomorphic
global sections $\p^0,\dots,\p^k$ of the line bundle $\F$
corresponding to the coordinate divisors of $F$. Taking into
account Proposition \ref{t-bund1}, we get conditions a) and b).
Further, fix $\O'\Subset\O$ and form a finite open covering
$\{\o_\a\}$ of $\hR\times\bar\O'$ such that for each $\o_\a$ at
least one of the sections $\p^l$ has no zeroes. Then choose an
open subcovering $\{\o_\a'\}$ of $\hR\times\bar\O'$ such that
$\o_\a'\Subset\o_\a$ for all $\a$. To prove c), we take $\d>0$
such that the mapping $j:\,\R^m\mapsto\hR$ maps every ball
$B(z',\d)$ with center $z'=x'+iy'\in T_{\O'},\ (jx',y')\in\o_\a'$
into the set $\o_\a$.

Conversely, let conditions a), b), c) be fulfilled. It follows
from Theorem \ref{t-bund} that global holomorphic sections $\p^0,\dots,\p^k$
of line bundles $\F_0,\dots,\F_k$ over $\hRO$ correspond to the divisors
$d_0,\dots,d_k$, respectively.
Note that if holomorphic sections $\{\p_\a\}$, $\{\tilde\p_\a\}$
of line bundles $\F$, $\tilde\F$ with transition functions
$\{g_{\a,\beta}\}$, $\{\tilde g_{\a,\beta}\}$ correspond to
divisors $d$, $\tilde d$ with a same Chern class, then the line
bundle with the transition functions $\{g_{\a,\beta}/\tilde
g_{\a,\beta}\}$ has a global holomorphic section $\{\p'_\a\}$
vanishing nowhere, therefore the sections $\{\p_\a\}$,
$\{\tilde\p_\a\p'_\a\}$  of the same line bundle $\F$
are assigned to the divisors $d$, $\tilde d$. Therefore we may
suppose that $\F_0=\F_1=\dots=\F_k=\F$.

Assume that all the sections $\p^l$
have a common zero $(\tau_0,y_0)$. Then in some neighborhood of
this point the sections $\p^l$ are defined by the functions
$\p^l_\a,\ l=0,\dots,k$. Choose a sequence $x_n\in\R^m$ such that
$jx_n\to\tau_0$ as $n\to\infty$. It is clear that the holomorphic
in $x+iy$ functions $\p^l_\a(j(x_n+x),y)$ converge as $n\to\infty$, uniformly
in some neighborhood of the point $(0,y_0)$, to the functions
$\p^l_\a(\tau_0+jx,y),\ l=0,\dots,k$, respectively. Hence Hurvitz'
Theorem implies that for arbitrary small $\d>0$ and sufficiently
large $n$ all the functions $\p^l_\a(jx,y)$ have zeroes in the
ball $B(x_n+iy_0,\d)$. This contradicts condition c). Thus all
the sections $\p^l$ have no common zeroes, and Proposition
\ref{t-map1} yields that there exists an almost periodic mapping with
the coordinate divisors $d_0,\dots,d_k$.\bs

\medskip

Now we prove a theorem which gives a method for constructing
holomorphic almost periodic mappings.

\begin{th}\label{t-prod}
\label{t-map3} Let $F(z)=[f^0(z):f^1(z):\dots:f^k(z)], \ \tilde
F(z)=[\tilde f^0(z):\tilde f^1(z):\dots:\tilde f^k(z)]$ be
holomorphic almost periodic mappings from $T_\O$ to $\CP$, a function $h\in
H(T_\O)$, and let the functions $f^l(z)\tilde f^l(z)/h(z),\
l=0,\dots,k$, be holomorphic on $T_\O$ and their divisors satisfy
condition c) of Theorem \ref{t-map2}. Then these functions are
homogeneous coordinates for some holomorphic almost periodic mapping.
\end{th}

\ws It follows from Proposition \ref{t-map1} that there exist
line bundles $\F$ and $\tilde\F$ over $\hRO$ with transition
functions $\{g_{\a\beta}\}$ and $\{\tilde g_{\a\beta}\}$, whose
sections $\{\p^l_\a\}$, $\{\tilde\p^l_\a\}$, $l=0,\dots,k$,
correspond to the coordinate divisors of the mappings $F$,
$\tilde F$. Then $\{\p^l_\a\tilde\p^l_\a\}$,
$l=0,\dots,k$, are global holomorphic sections of the line bundle
with the transition functions $\{g_{\a\beta}\tilde
g_{\a\beta}\}$. Fix a point $(\tau_0,y_0)\in\hRO$ and its
neighborhood $\o$ such that the sections
$\{\p^l_\a\tilde\p^l_\a\}$ are defined by holomorphic functions
$q^l(\tau,y)$, $l=0,\dots,k$, on $\o$. Choose a sequence
$x_n\in\R^m$ such that $jx_n\to\tau_0$ as $n\to\infty$. There
exists $\d>0$ such that for $n$ sufficiently large the image of
every ball $B(x_n+iy_0,\d)$ under the map $j$ lies in $\o$ and at
least one of the holomorphic functions $q^l(jx,y)/h(x+iy)$ has no
zeroes in each ball. For the sake of being definite, suppose that
for an infinite sequence of numbers $n$ the functions $q^0(jx,y)$
and $h(x+iy)$ have the same zeroes in $B(x_n+iy_0,\d)$. Since
$h(x+iy)$ divides all the functions $q^l(jx,y)$, so does the
function $q^0(jx,y)$. Using Proposition \ref{t-quot}, we get that
in a sufficiently small neighborhood $\o'$ of the point
$(\tau_0,y_0)$,
$$
q^l(\tau,y)=q^0(\tau,y)r^l(\tau,y),\quad l=1,\dots,k,
$$
with some functions $r^l(\tau,y)$ holomorphic in $\o'$.

Choose an open covering $\{\o_\a\}$ of the set $\hRO$ such that
for each $\a$ there exists a number $l(\a)$ and functions
$q^l_\a$, $r^l_\a$, holomorphic in $\o_\a$  such that $q^l_\a$
define the sections $\{\p^l_\a\tilde\p^l_\a\}$, $l=0,\dots,k$.
We have
$$
q^l_\a(\tau,y)=q^{l(\a)}_\a(\tau,y)r^l_\a(\tau,y),\quad l=0,\dots,k.
$$
It is clear that the functions $q^{l(\beta)}_\beta/q^{l(\a)}_\a$
are holomorphic and have no zeroes on $\o_\a\cap\o_\beta$, therefore
the functions
$\{g_{\a\beta}\tilde g_{\a\beta}q^{l(\beta)}_\beta/q^{l(\a)}_\a\}$
are transition functions for some line bundle $\F'$.
Then the functions $\{r^l_\a\}$ form $k+1$ holomorphic sections of
this bundle and $r_\a^{l(\a)}\equiv 1$ on $\o_\a$. Using Proposition \ref{t-map1},
we obtain that the line bundle $\F'$ generates a holomorphic
almost periodic mapping $H$ with the homogeneous coordinates
$f^l(z)\tilde f^l(z)/h(z),\ l=0,\dots,k$. \bs

\begin{cor}
\label{t-map4} (for $m=1$, see \cite{FP}). The product of two
meromorphic almost periodic  functions is almost periodic if and
only if the zeroes and poles of the product are uniformly separated
in every $T_{\O'}$ with $\O'\Subset\O$.
\end{cor}

Note that there exists a  holomorphic almost periodic mapping $F$
such that for every homogeneous holomorphic representation
$F(z)=[f^0(z):f^1(z):\dots:f^k(z)]$, none of the functions $f^l$
belongs to $APH(T_\O)$.

Indeed,  let $\phi(\zeta)$ be a holomorphic function on the
plane with zeroes at all points with integer coordinates. Fix
points $\zeta_1,\dots,\zeta_k$ in the plane such that neither the
points, nor their differences have integer coordinates. Then the
divisors $d_0,\dots,d_k$ of the functions
$h_0(z^1,z^2)=\phi(z^1-iz^2)$, $h_1(z^1,z^2)=\phi(z^1-iz^2+\zeta_1),
\dots,h_k(z^1,z^2)=\phi(z^1-iz^2+\zeta_k)$ are periodic, have
the same nonzero Chern class, and satisfy the conditions c) of
Theorem \ref{t-map2}. Therefore there
exists a holomorphic almost periodic mapping with the coordinate
divisors $d_0,\dots,d_k$, but there exist no functions from
$APH(T_\O)$ with these divisors.

Nevertheless, the following theorem is true.

\begin{th} Let $F(z)$ be a holomorphic almost periodic mapping of
$T_\O$ into the projective space $\CP$. Then there exists an
almost periodic divisor $d$ in $T_\O$ and $g^0(z),\dots,g^k(z)\in
APH(T_\O)$ such that their common zeros are contained in the
support of $d$, and \be\label{rep} F(z)=[g^0(z):\dots:g^k(z)] \en
on $T_\O\setminus\supp \, d$; the representation (\ref{rep}) with
holomorphic almost periodic functions without common zeroes exists
if and only if the coordinate divisors of $F$ have zero Chern
class.
\end{th}

\ws If all the functions $g^j$ in (\ref{rep}) have no common
zeroes, then their divisors are the coordinate divisors of $F$.
It follows from Theorem \ref{t-Cher} that all these divisors
have zero Chern class.

Further, let $[f^0(z):\dots:f^k(z)]$ be a homogeneous
representation of $F$ with functions $f^j$ holomorphic in $T_\O$
without common zeroes. Let $d_0$ be the divisor of $f^0(z)$. Using
Theorem \ref{t-comp}, take an almost periodic divisor $d$ such
that $d_0+d$ is the divisor of some $g^0\in APH(T_\O)$; we take
$d=0$ in the case $c(d_0)=0$. Since the mapping $\tilde
F(z)=[1:g^0(z):\dots:g^0(z)]$ is almost periodic, Theorem
\ref{t-prod} with $h(z)=f^0(z)$ implies that the mapping
$[1:g^1(z):\dots:g^k(z)]$ with $g^j(z)=f^j(z)g^0(z)/f^0(z),\
j=1,\dots,k$, is almost periodic, too. Hence the functions
$g^j(z),\ j=0,1,\dots,k$, are almost periodic and have common
zeroes only in $\supp d$. It is clear that the representation
(\ref{rep}) is valid on the set $T_\O\setminus\supp \,d$. \bs

\bigskip
\section{Further remarks}
\label{sec7}

Our method of determining of holomorphic function by almost
periodic divisor does not allow to control the growth of the
function (even in the case $T_\O=\C$). Perhaps, we need a more
precise integral representation for the $\bar\partial$-problem.

Note also that it looks natural to replace the group of mappings
$\{z\mapsto z+t\}_{t\in\R^m}$ in the definition of almost periodic
functions, to any group of automorphisms; apparently, then we will
need another compactification instead of Bohr's one and another
integral representation for the $\bar\partial$-problem.

\bigskip

{\bf Acknowledgements.} The author is grateful to Professor
G.\,Henkin and Professor A.\,Rashkovskii for useful discussions.

\bigskip
{\it Department of Function Theory and Functional Analysis

Kharkiv National University

Svobody sq, 4, Kharkiv 61077, Ukraine

\medskip
e-mail: favorov@assa.vl.net.ua}

\end{document}